\documentclass[a4 paper, 12pt]{amsart}

\usepackage{amsmath,amssymb,amscd,amsfonts,verbatim}

\usepackage[mathscr]{eucal}

\newcommand{\pionealg}{\pi_1^{\mathrm{alg}}}
\newtheorem{thm}{Theorem}[section]
\newtheorem{lem}[thm]{Lemma}
\newtheorem{rem}[thm]{Remark}
\newtheorem{prop}[thm]{Proposition}
\newtheorem{cor}[thm]{Corollary}

\theoremstyle{remark}

\newcommand{\iso}{\cong}
\newcommand{\into}{\hookrightarrow}

\newcommand{\Span}[1]{\left<#1\right>}
\newcommand{\C}{\mathbb C}
\newcommand{\Z}{\mathbb Z}
\newcommand{\Q}{\mathbb Q}

\newcommand{\F}{\mathbb F}

\newcommand{\sM}{\mathcal M}
\newcommand{\pp}{\mathbb P}
\DeclareMathOperator{\Aut}{Aut}
\DeclareMathOperator{\Pic}{Pic}

\DeclareMathOperator{\Tors}{Tors}

\DeclareMathOperator{\Proj}{Proj}
\newcommand{\epsi}{\varepsilon}

\newcommand{\ga}{\gamma}

\newcommand{\om}{\omega}
\newcommand{\Ga}{\Gamma}
\newcommand{\De}{\Delta}
\newcommand{\Si}{\Sigma}
\newcommand{\si}{\sigma}
\newcommand{\ze}{\zeta}
\newcommand{\fie}{\varphi}

\newcommand{\OO}{\mathcal{O}}

\newcommand{\inv}{^{-1}}

\numberwithin{equation}{section}
\title[Numerical Campedelli surfaces\dots ]{Numerical Campedelli
surfaces with fundamental group of order 9}
\author{Margarida Mendes Lopes and Rita Pardini}
\date{}
\thanks{The first author is a member of the Center for Mathematical
Analysis, Geometry and Dynamical Systems and the second author is
a member of G.N.S.A.G.A.-I.N.d.A.M. This research was partially
supported by the Italian project ``Geometria sulle variet\`a
algebriche" (PRIN COFIN 2004) and by FCT (Portugal) through
program POCTI/FEDER and Project POCTI/MAT/44068/2002.}

\begin{document}

\begin{abstract}
We give explicit constructions of all the numerical Cam\-pe\-delli
surfaces, i.e.\ the minimal surfaces of general type with $K^2=2$
and $p_g=0$, whose fundamental group has order 9. There are three
families, one with $\pionealg=\Z_9$ and two with
$\pionealg=\Z_3^2$.

We also determine the base locus of the bicanonical system of
these surfaces. It turns out that for the surfaces with
$\pionealg=\Z_9$ and for one of the families of surfaces with
$\pionealg=\Z_3^2$ the base locus consists of two points. To our
knowlegde, these are the only known examples of surfaces of
general type with $K^2>1$ whose bicanonical system has base
points.

\end{abstract}
\maketitle

\section{Introduction}

A numerical Campedelli surface is a minimal surface of general
type with $K^2=2$ and $p_g=0$. It is known (cf.\
\cite[Chap.~VII.10]{bpv2}) that the algebraic fundamental group
$\pionealg$ of such a surface is finite, of order at most 9. In
this paper, we construct explicitly all the numerical Campedelli
surfaces with $|\pionealg|=9$.

For this, given a surface $S$ as above, we consider its universal
cover $Y\to S$. Assume for simplicity that $K_S$ is ample. Then by
results of Konno (\cite{konnopg}) the canonical map of $Y$ is an
isomorphism and $Y$ is the intersection of a normal threefold
$W\subset \pp^7$ of degree 5 or 6 with a cubic hypersurface. Using
Fujita's classification of these threefolds
(\cite{fujita1,fujita2,fujita3}), we prove that in our situation
$W$ actually has degree 6 and that there are three possibilities
for $W$, each giving rise to a family of numerical Campedelli
surfaces with algebraic fundamental group of order 9
(Theorem~\ref{main}).

More precisely (cf.\ \S3), there is an irreducible family of
dimension 6 of numerical Campedelli surfaces with
$\pionealg=\Z_9$, and two irreducible families, of dimension
respectively 7 and 6, of numerical Campedelli surfaces with
$\pionealg=\Z_3^2$. Recall that the expected number of moduli of a
numerical Campedelli surface is 6; as far as we know, these are
the only known examples in which the dimension of the moduli space
is greater that the expected one.

Furthermore we prove (Theorem~\ref{moduli}) that the moduli space
of these surfaces has exactly two connected components, both
irreducible, corresponding to $\pionealg=\Z_9$ and
$\pionealg=\Z_3^2$. We also show that for all these surfaces
$\pionealg$ coincides with the topological fundamental group
(Theorem~\ref{pi1}).

We determine also the base locus $\Ga$ of the bicanonical
system $|2K_S|$ (Theorem~\ref{basepoints}). If
$\pionealg(S)=\Z_9$, then $\Ga$ consists of two points. If
$\pionealg(S)=\Z_3^2$, then $\Ga$ is empty in general but there is
a codimension~1 subvariety of the moduli space where $\Ga$ consist
of two points. This fact is quite surprising. Indeed, excluding
the case $K^2=1$, $p_g=0$ in which $|2K_S|$ is a pencil, the
bicanonical map of a minimal surface of general type is
generically finite (\cite{xiaocan}) and, by the work of various
authors (cf.\ \cite{ciro}), it is known to be a morphism if either
$p_g>0$ or $K^2>4$. To our knowledge, the examples presented here
are the only known surfaces of general type with $K^2>1$ whose
bicanonical map is not a morphism.

The paper is organized as follows. In \S2 we establish some preliminaries. In \S3 we state the classification  theorem (whose proof will be given in \S5) and describe in detail the three  families of Campedelli surfaces with $\pionealg$
 of order 9. In \S4 we study the moduli and the bicanonical system
of such surfaces and finally in \S5 we prove the classification theorem.

\medskip

\noindent {\bf Notation} We work over the complex numbers. All
varieties are projective algebraic. All the notation we use is
standard in algebraic geometry. We just recall the definition of
the numerical invariants of a smooth surface $S$: the
self-intersection number $K^2_S$ of the canonical divisor $K_S$,
the {\em geometric genus} $p_g(S):=h^0(K_S)=h^2(\OO_S)$, the {\em
irregularity} $q(S):=h^0(\Omega^1_S)=h^1(\OO_S)$ and the {\em
holomorphic Euler characteristic} $\chi(S):=1+p_g(S)-q(S)$.

\medskip

\noindent {\bf Acknowledgments} We are very grateful  to Miles Reid for all his suggestions (mathematical, linguistic, typographical), which improved substantially the presentation of this paper. In particular we wish to thank him for pointing out the unified presentation of surfaces of type B1 and B2 (see \ref{common}). 

We are also indebted to the referee for his thorough reading of the paper and for suggesting further improvements of the presentation.

\section{Set-up and preliminaries}

Here we establish the notation and recall or prove some facts
that we use throughout all the paper.

Let $S$ denote a numerical Campedelli surface, i.e.\ a minimal
projective surface of general type with $K^2_S=2$ and $p_g(S)=0$
(and so $q(S)=0$) and let $G$ denote the algebraic fundamental
group of $S$.
\begin{prop}[\cite{milesk2}, cf.\ \cite{3chi}] \label{pi9}
The group $G$ is finite, of order $\le 9$.
\end{prop}

We assume from now on that $G$ has order 9 and we consider the
corresponding \'etale $G$-cover $\pi\colon Y\to S$.
 \begin{prop}\label{Y}
 \begin{enumerate}
\item $Y$ is a smooth minimal surface of general type with
$K^2_Y\!=\!18$, $p_g(Y)\!=\!8$, $q(Y)\!=\!0$, hence
$K_Y^2\!=\!3p_g(Y)-6$;

\item the canonical map $\fie\colon Y\to\pp^7$ is birational onto
its image $V$.
 \end{enumerate}
 \end{prop}
 \begin{proof}
Since $\pi$ is \'etale, $Y$ is smooth and $K_Y=\pi^*K_S$ is nef
and big. Hence $Y$ is minimal of general type and
$K^2_Y=9K^2_S=18$, $\chi(Y)=9\chi(S)=9$. By Proposition~\ref{pi9},
we have $\pionealg(Y)=0$ and thus, in particular, $q(Y)=0$. It
follows that $p_g(Y)=8$ and so $K_Y^2=3p_g(Y)-6$.

Suppose that the canonical map $\fie$ of $Y$ is not birational.
Then, by \cite{konnopg}, $\fie$ is generically finite of degree
$2$ and its image is a ruled surface, which is rational by the
regularity of $Y$. But then, by \cite[Cor. 5.8]{beauville} (cf.\
\cite[Proposition~4.1]{3chi}), $G=\Z_2^r$, contradicting the
assumption that $G$ has order 9.
 \end{proof}

Conversely one has:

\begin{prop}\label{campe}
Let $V$ be a surface with canonical singularities and ample $K_V$,
and with $K^2_V=18$, $p_g(V)=8$, $q(V)=0$; let $G$ be a group of
order 9 that acts freely on $V$.

Then the quotient surface $S:=V/G$ is the canonical model of a
numerical Campedelli surface with $\pionealg=G$.
 \end{prop}
\begin{proof} Since the quotient map $p\colon V\to S$ is \'etale,
$S$ has canonical singularities and $K_V=p^*K_S$. Since $K_V$ is
ample, $K_S$ is also ample, and therefore $S$ is the canonical
model of a surface of general type. Since $p$ is \'etale, we have
$K^2_S=K^2_V/9=2$, $\chi(S)=\chi(V)/9=1$. In addition, $q(S)=0$
since $q(V)=0$, and thus $p_g(S)=0$. Let $S'\to S$ and $Y\to V$ be
the minimal resolutions. The surface $S'$ is a numerical
Campedelli surface and the map $p\colon V\to S$ induces an \'etale
map $p'\colon Y\to S'$ with Galois group $G$. Therefore there is a
surjection $\pionealg (S')\twoheadrightarrow G$. By
Proposition~\ref{pi9} this is actually an isomorphism.
\end{proof}

The next result  is needed in \S3 to explain the choices of  linearization in the construction of the examples. It is well known to experts (cf.\ \cite{milesk2}),
but we recall it here for lack of a published reference. 

 \begin{prop}\label{Grep}The representation of $G$ on
$H^0(Y, K_Y)$ decomposes as the direct sum of the~8 nontrivial
characters of $G$.
 \end{prop}

 \begin{proof}
Let $\chi\in G^*$ be a character of $G$ and $V_{\chi}\subset
H^0(Y,K_Y)$ the corresponding eigenspace. There is a natural
isomorphism $G^*\to\Tors S$, the torsion subgroup of $\Pic S$,
which induces a natural identification $V_{\chi}\iso H^0(S,
K_S+\eta_{\chi})$, where $\eta_{\chi}$ is the line bundle
associated to $\chi\in G^*$. Notice that $V_1=\{0\}$, since we
have $p_g(S)=0$. By Proposition~\ref{pi9}, there is no irregular
\'etale cover of $S$, hence for every $\eta\in \Tors S$ we have
$h^1(S, \eta)=0$. Since $\chi(S)=1$, Serre duality gives $h^0(S,
K_S+\eta)=1$ for every $\eta\in \Tors S\setminus \{0\}$, i.e.\
$\dim V_{\chi}=1$ for every $\chi\in G^*\setminus\{0\}$.
\end{proof}

 \section{The classification}

 The main result of this paper, which will be proven in \S\ref{secproof}, is the following classification theorem:
\begin{thm}\label{main}
Let $S$ be a numerical Campedelli surface with $\pionealg(S)$ of
order~9. Then:
\begin{enumerate}
\item if $\pionealg(S)\iso \Z_9$, then $S$ is a surface of type A
(cf.\ \S\ref{ssecA}); \item if $\pionealg(S)\iso \Z_3^2$, then $S$ is a
surface of type B1 or B2 (cf.\ \S\ref{ssecB1}, \S\ref{ssecB2}).
\end{enumerate}
\end{thm}

 In this section we describe in detail  the three families of Campedelli surfaces with $\pionealg$ of order 9 (i.e., type A, type B1, type B2). 
 The notation is consistent with the previous
section. 
\smallskip

The examples are obtained as  quotients of
surfaces $Y$ with $K^2_Y=18$, $p_g(Y)=8$ and $q(Y)=0$ by a group
$G$ of order 9 acting freely.

In each case $Y$ is the minimal resolution of a surface $V$  of $\pp^7$ with canonical singularities, obtained as the intersection of a normal  threefold $W\subset \pp^7$ of degree 6 with a cubic hypersurface. The  line bundle $\OO_W(1)$ is equal to $-2K_W$ and the surface $V$ is contained in the smooth part of $W$, hence by adjunction  $K_V=\OO_V(1)$. ($W$ is actually smooth in cases A and B1 and it is a cone over a smooth surface of $\pp^6$ in case B2). 
 It is then easy to check that $V$ has the right invariants:  $K^2_V=18$, $p_g(V)=8$, $q(V)=0$. By construction, the threefold $W$ is the intersection of all the quadrics containing $V$.

The group $G$ acts on $W\subset \pp^7$ and the action is free on $W$ outside a finite set.  We choose   a linearization of $\OO_W(3)$ and we let $T_1\subset H^0(\OO_W(3))$ be the space of invariant sections with respect to this linearization.  Then we show that the linear system $|T_1|$ is free on $W$, so that the general $V\in |T_1|$ satisfies the assumptions of Proposition \ref{campe} and $V/G$ is the canonical model of a Campedelli surface with $\pionealg=G$.
  In all the examples, the map $H^0(\OO_{\pp^7}(3))\to H^0(\OO_W(3))$ is surjective, so that $V$ is cut out on $W$ by a cubic hypersurface of $\pp^7$.

\begin{rem}{\em 
In the above construction, the chosen  linearization of $\OO_W(3)$ induces a decomposition into $G$-eigenspaces: $$H^0(\OO_W(3))=\oplus_{\chi\in G^*} T_{\chi}.$$ For every $\chi\in G^*$ the linear system $|T_{\chi}|$ consists of $G$-invariant surfaces, hence one could  hope to obtain surfaces with the required properties also for some character $\chi\ne 1$.  It turns out (cf. Lemma \ref{system}) that  this is not possible  because for every $\chi\ne 1$ the system $|T_{\chi}|$ has base points  with non trivial stabilizer in $G$, and therefore $G$ does not act freely on the surfaces of $|T_{\chi}|$ for $\chi\ne 1$. }
\end{rem}
\begin{rem} {\em  In  each example we compute  the dimension of the invariant subspace  $T_1$ by writing down a basis of it.  Alternatively, the dimension of $T_1$ can be computed by using   the equivariant Riemann-Roch theorem (\cite{ypg}, Corollary in (8.6)).  }
\end{rem}

\subsection{Surfaces of type A}\label{ssecA}In this example $G\iso\Z_9$.

We take $W:=\pp^1\times\pp^1\times \pp^1$ and consider homogeneous
coordinates $(x_0, x_1),(y_0,y_1),(z_0,z_1)$ on $W$ and the
corresponding affine coordinates $x:=x_1/x_0$, $y:=y_1/y_0$ and
$z:=z_1/z_0$.

We fix a primitive third root $\om$ of 1 and we let a generator
$g$ of $G$ act  by
\begin{equation}\label{Z9action}
g\colon(x,y,z)\mapsto(y,z,\om x).
\end{equation}
Note that the fixed points of $g^3$ are the 8 points with the
three affine coordinates $x,y,z$ equal either to $0$ or to
$\infty$.

Set $H:=\OO_W(1, 1,1)$. The linear system $|H|$ embeds $W$ in
$\pp^7$ as a smooth threefold of degree 6.
The action of $G$ extends to $\pp^7$ and the projective representation of $G$ on $\pp^7$ induces a linear representation on $H^0(\OO_{\pp^7}(3))$, hence a linearization of $3H$.
Let $T_1\subset H^0(W, 3H)$ be the subspace of $G$-invariant elements. In affine coordinates  a basis of $T_1$ is given by:
\begin{gather}\label{basisA}
1, \ x^3+y^3+z^3, \ x^2y+y^2z+\omega z^2x, \ x^2z+\omega y^2 x+\omega z^2y\\
x^3y^3+y^3z^3+z^3x^3,\  x^3y^2z+\omega y^3z^2x+z^3x^2y,\nonumber
\\ x^3yz^2+\omega^2 y^3zx^2+z^3xy^2, \ x^3y^3z^3.\nonumber
\end{gather}
Using this basis it is easy to check that the system $|T_1|$ is free of dimension 7. 

 By Bertini's theorem, the general element of
$|T_1|$ is a smooth surface on which $G$ acts freely. Let
$V\in|T_1|$ be a surface with at most rational double points and
not passing through the fixed points of $g^3$. By adjunction, the
canonical divisor $K_V$ is the restriction of $H$ to $V$, hence
$K^2_V=3H^3=18$. The adjunction sequence gives: $p_g(V)=8$,
$\chi(V)=9$, $q(V)=0$. Set $S:=V/G$. Then by
Proposition~\ref{campe} the surface  $S$ is the canonical model of a numerical
Campedelli surface with $\pionealg =G$.

\begin{prop} \label{moduliA}The family of surfaces of type $A$
depends on~$6$ moduli.
\end {prop}

\begin{proof} Consider $V_1,V_2\in|T_1|$  such that the surfaces $S_1:=V_1/G$ and $S_2:=V_2/G$ are isomorphic. Since $V_i$ is the (algebraic) universal cover of $S_i$,  the isomorphism $S_1\to S_2$ lifts to an isomorphism $V_1\to V_2$.  Since $V_1$ and $V_2$ are canonically embedded in $\pp^7$,  the isomorphism $V_1\to V_2$ is induced by an automorphism $\ga$ of $\pp^7$. The threefold $W\subset \pp^7$ is the intersection of all the quadrics containing $V_i$, $i=1,2$ (cf. \cite[Theorem 3.1]{konnopg}), and therefore $\ga$ maps $W$ to itself.  Since every automorphism of $W$ extends to $\pp^7$, we may regard $\ga$ as an automorphism of $W$. 

By construction, $\ga$ belongs to the normalizer $\Ga$ of $G$ in $\Aut W $. Since the group $G$ is finite, the connected  component $\Ga_1$
  of the identity
in $\Ga$ is actually contained in the centralizer of $G$. In addition, if   $\ga$ is an element of  $\Ga_1$, then by continuity $\ga$ does not permute the three copies of $\pp^1$.  Then it is easy to see that $\Ga_1$ consists of the maps of the form $(x,y,z)\mapsto (\lambda
x, \lambda y, \lambda z)$ for $\lambda\in \C^*$. Hence we obtain a
family of Campedelli surfaces with 6 moduli.
\end{proof}
\begin{rem}\label{expdim}{\em 
The expected dimension of the moduli space of a minimal surface of general type $S$ is $10\chi(S)-2K^2_S=\chi(T_S)$. In particular  the expected dimension of the moduli space
of Campedelli surfaces is equal to 6.}
\end{rem}
\begin{rem} \rm
J.H. Keum has kindly communicated to us an example, due to Persson, of
a numerical Campedelli surface with $\pionealg\iso\Z_9$. One can
check that this example is given by the above construction.
\end{rem}

\subsection{Surfaces of type B1}\label{ssecB1} In this example
$G\iso\Z_3^2$.

We take $W$ to be the  flag variety $\{x_0y_0+x_1y_1+x_2y_2=0\}\subset
\pp^2\times {\pp^2}^*$.  We denote by $H$ the restriction to $W$ of $\OO_{\pp^2\times{\pp^2}^*}(1,1)$, so that $W$ is embedded by $|H|$ into $\pp^7$ 
as a smooth threefold of degree 6.
Let $g_1,g_2\in G$ be generators acting on $\pp^2$ by
\begin{equation}\label{Z3action}
(x_0,x_1,x_2)\overset{g_1}{\mapsto}(x_0,\om x_1, \om^2x_2), \quad
(x_0,x_1,x_2)\overset{g_2}{\mapsto}(x_1,x_2,x_0),
\end{equation}
 where $\om\ne 1$ is a third root of 1.
The induced $G$-action on $\pp^2\times {\pp^2}^*$ obviously
preserves $W$ and it is easy to  check that every  nontrivial element of $G$ fixes finitely many points of $W$.

Notice that, although the action of $G$ on $\pp^2$ is not induced
by a linear action, the corresponding action on $\pp^2\times
{\pp^2}^*$ is induced by a linear $G$-action on
$H^0(\OO_{\pp^2\times{\pp^2}^*}(1,1))$.
In particular, $G$ acts on $\OO_{\pp^2\times{\pp^2}^*}(1,1)$, on  $H$ and on their  multiples. 
There is an exact sequence:
$$0\to H^0(\OO_{\pp^2\times{\pp^2}^*}(2,2))\overset{i}{\to} H^0(\OO_{\pp^2\times{\pp^2}^*}(3,3))\to H^0(W, 3H)\to 0.$$
The map $i$ is multiplication by the invariant section $x_0y_0+x_1y_1+x_2y_2$ and thus it is equivariant with respect to the action of $G$.  In particular,  the subspace $T_1\subset H^0(W,3H)$ on which $G$ acts trivially is the image of the subspace $R_1\subset H^0(\OO_{\pp^2\times\pp^2}(3,3))$ on which $G$ acts trivially. A basis of $R_1$, which has dimension 12,  is the following:
\begin{gather*}
(x_0^3+x_1^3+x_2^3)(y_0^3+y_1^3+y_2^3), (x_0^3+x_1^3+x_2^3)(y_0y_1y_2),\\
(x_0x_1x_2)(y_0^3+y_1^3+y_2^3), (x_0x_1x_2)(y_0y_1y_2), \\
(x_0^2x_1+x_1^2x_2+x_2^2x_0)(y_0^2y_1+y_1^2y_2+y_2^2y_0),\\
(x_1^2x_0+x_2^2x_1+x_0^2x_2)(y_1^2y_0+y_2^2y_1+y_0^2y_2),\\
(x_0^2x_1+ \omega x_1^2x_2+\omega^2x_2^2x_0)( y_0^2y_1+ \omega^2 y_1^2y_2+\omega y_2^2y_0),\\
( x_1^2x_0+\omega x_2^2x_1+\omega^2x_0^2x_2)(y_1^2y_0+\omega^2y_2^2y_1+\omega y_0^2y_2),\\
(x_0^2x_1+ \omega^2 x_1^2x_2+\omega x_2^2x_0)( y_0^2y_1+ \omega y_1^2y_2+\omega^2 y_2^2y_0),\\
( x_1^2x_0+\omega^2 x_2^2x_1+\omega x_0^2x_2)(y_1^2y_0+\omega y_2^2y_1+\omega^2 y_0^2y_2),\\
(x_0^3+\omega x_1^3+\omega^2x_2^3)(y_0^3+\omega^2y_1^3+\omega y_2^3), (x_0^3+\omega^2x_1^3+\omega x_2^3)(y_0^3+\omega y_1^3+\omega^2y_2^3).
\end{gather*}
Next we write down a basis of the  subspace of $H^0(\OO_{\pp^2\times{\pp^2}^*}(2,2))$ on which $G$ acts trivially:
\begin{gather*}
x_0^2y_0^2+x_1^2y_1^2+x_2^2y_2^2,\\ x_0^2y_1y_2+x_1^2y_0y_2+x_2^2y_0y_1,  x_1x_2y_0^2+x_0x_2y_1^2+x_0x_1y_2^2,\\
x_0x_1y_0y_1+x_1x_2y_1y_2+x_2x_0y_2y_0.
\end{gather*}
It follows that the space $T_1$ has dimension $12-4= 8$. Let $z_0, \dots z_8$ be a basis of $H^0(\OO_{\pp^2\times{\pp^2}^*}(1,1))$ on which $G$ acts diagonally. Then $z_0^3, \dots z_8^3$ restrict on $W$ to elements of  $T_1$, hence the linear system 
$|T_1|$ is free.

By  Bertini's theorem,
the  general $V\in |T_1|$ is a smooth surface on which $G$ acts freely.
Let $V\in |T_1|$ be a surface with at most rational double points
and not passing through any point fixed by a nontrivial element of
$G$. By adjunction, the canonical divisor $K_V$ is the restriction
of $H$ to $V$, hence $K^2_V=3H^3 =18$. The adjunction sequence for
$V\subset W$ gives: $p_g(V)=8$, $\chi(V)=9$, $q(V)=0$. Set
$S:=V/G$. Then by Proposition~\ref{campe} the surface $S$ is the
canonical model of a numerical Campedelli surface with
$\pionealg=G$.

To compute the number of moduli for surfaces of type B1 we need the following description of $\Aut W$:
\begin{lem}\label{autW}
There is an exact sequence:
\[
0\to\Aut\pp^2\stackrel{i}{\longrightarrow}\Aut W
\stackrel{\si}{\longrightarrow}\Z_2\to 0.
\]
In particular $G$, being of order 9, can be identified with a
subgroup of $\Aut\pp^2$.
\end{lem}
\begin{proof}
Every automorphism $g$ of $\pp^2$ extends naturally to an automorphism of $\pp^2\times {\pp^2}^*$ that preserves the flag variety $W$.
The Picard group of $W$ is isomorphic to $\Z^2$, since $W$ is a
$\pp^1$-bundle over $\pp^2$. Set $H_1:=\OO_W(1,0)$,
$H_2:=\OO_W(0,1)$. The classes of $H_1$ and $H_2$ are a basis of
$H^2(W,\Q)$ and they satisfy $H_1^3=H_2^3=0$,
$H_1K_W^2=H_2K_W^2=12$. It is easy to check that these numerical
conditions characterize $H_1$ and $H_2$. Hence every automorphism
$\ga$ of $W$ either preserves the classes of $H_1$ and $H_2$ or it
exchanges them. The map $\si\colon \Aut W\to\Z_2$ sends $\ga$ to
$0$ in the former case and to $1$ in the latter case. Clearly the
map $\si$ is not trivial and its kernel contains $\Aut\pp^2$.

Now let $g\in\ker\si$. Then $g$ induces an automorphism $g_1$ of
$\pp^2=|H_1|^*$ and an automorphism $g_2$ of ${\pp^2}^*=|H_2|^*$
such that the automorphism $g_1\times g_2$ of $\pp^2\times
{\pp^2}^*$ restricts to $g$ on $W$. Since the flag variety
$W$ is preserved, it follows that $g_2$ is equal to $^t g_1\inv$,  hence $g$ is in the subgroup
$\Aut\pp^2$.
\end{proof}

\begin{prop}\label{moduliB1} The family of surfaces of type B1
depends on $7$ moduli.
\end{prop}

\begin{proof}
By the argument in the proof of Proposition \ref{moduliA}, it is enough to
 show that the  centralizer  $\Ga$ of $G$ in $\Aut W$  is discrete.

Let $\Ga_1$ be the connected component of the identity in $\Ga$ and let $\ga\in \Ga_1$.
By Lemma~\ref{autW} and by continuity, $\ga$ is induced by an
element of $\Aut\pp^2$, that we again denote by $\ga$. Since $\ga g\ga\inv=g$  for every $g\in G$, $\ga$ permutes the fixed points of $g$ hence, again by continuity, the fixed points of $g$ on $\pp^2$ are also fixed points of $\ga$.
Since the set of points in $\pp^2$ that are fixed by a nontrivial
$g\in G$ consists of 12 points in general position, $\ga$ is the
identity.
\end{proof}

\begin{rem}\rm As we recalled in Remark \ref{expdim}, the expected dimension of the moduli space of Cam\-pe\-delli surfaces is 6. To the best of our knowledge, the surfaces of type B1 are the only known example where the actual dimension exceeds the expected one.
\end{rem}

\begin{rem}\rm The first known example of a numerical Campedelli
surface with torsion $\Z_3^2$ is due to Xiao Gang,
\cite[Example~4.11]{xiao}. A particular case of the construction
above has been kindly communicated to us by JongHae Keum, who
attributes it to Xiao Gang and Beauville. See also \cite{cmp},
Example 2 of \S5.
\end{rem}

\subsection{Surfaces of type B2}\label{ssecB2} This is the second example with
$G\iso\Z_3^2$. It is also the only example in which $W$ is singular.

 Let $\epsi\colon \Si\to \pp^2$ be the blowup at the
three coordinate points. The $G$-action on $\pp^2$ defined 
 in \ref{Z3action} extends to a $G$-action on
$\Si$.

The action on $\pp^2$ does not lift to a linear action on
$H^0(\OO_{\pp^2}(1))$, but it induces a natural linear action of
$G$ on $H^0(\OO_{\pp^2}(3))$ that preserves the subspace $T$ of
sections vanishing at the coordinate points. Let $\chi_1\in G^*$
be the character such that $\chi_1(g_1)=\om$, $\chi_1(g_2)=1$ and
let $\chi_2$ be the character such that $\chi_2(g_1)=1$,
$\chi_2(g_2)=\om$. The following is a basis of the vector space
$T$ consisting of eigenvectors for the $G$-action. We write the
corresponding character next to each element:
\begin{align}\label{base}
 x_0^2x_1+x_1^2x_2+x_2^2x_0\quad &(1, 0);
 \quad x_1^2x_0+x_2^2x_1+x_0^2x_2\quad &(2,0); \\
 x_0^2x_1+ \om x_1^2x_2+\om^2x_2^2x_0 \quad &(1, 2);
 \quad x_1^2x_0+\om x_2^2x_1+\om^2x_0^2x_2 \quad &(2,2);
 \nonumber\\
 x_0^2x_1+ \om^2 x_1^2x_2+\om x_2^2x_0 \quad &(1, 1);
 \quad x_1^2x_0+\om^2x_2^2x_1+\om x_0^2x_2\quad &(2,1);
 \nonumber
 \\x_0x_1x_2\quad &(0,0).\nonumber
\end{align}
We fix homogeneous coordinates $(z_{ij})$ on $\pp^7$, for
$(i,j)\in \Z_3^2\setminus\{(0,0)\}$, and we let $G$ act on
$z_{ij}$ via the character $\chi_1^i\chi_2^j$. Let $P$ be the
point whose only nonzero coordinate is $z_{01}$ and identify
$\pp^6$ with the $G$-invariant hyperplane $z_{01}=0$. The rational
map $\pp^2\to \pp^6\subset \pp^7$ defined by
\begin{align}\label{map}
z_{12}= x_0^2x_1+x_1^2x_2+x_2^2x_0, &\quad z_{22}=
x_1^2x_0+x_2^2x_1+x_0^2x_2, \\
z_{11}= x_0^2x_1+ \om x_1^2x_2+\om^2x_2^2x_0, &\quad z_{21}=
x_1^2x_0+\om x_2^2x_1+\om^2x_0^2x_2, \nonumber \\
z_{10}= x_0^2x_1+ \om^2x_1^2x_2+\om x_2^2x_0, &\quad
z_{20}=x_1^2x_0+\om^2x_2^2x_1+\om x_0^2x_2, \nonumber \\
z_{02}=x_0x_1x_2 \nonumber
\end{align}
induces a $G$-equivariant embedding $\Si\to\pp^6\subset \pp^7$.
Notice that the embedding is not induced by a $G$-equivariant embedding $T\to H^0(\OO_{\pp^7}(1))$. Indeed by
(\ref{base}) the representation of $G$ on $T$ differs from the one induced by the
representation on $H^0(\OO_{\pp^7}(1))$ by multiplication by
$\chi_2^2$. This is due to the fact  that we have chosen the linearization of  $H^0(\OO_{\pp^7}(1))$ so that it will agree with the one  induced by the identification with $H^0(V, K_V)$,  where $V$ is any of the surfaces that we are constructing. Since by Proposition \ref{Grep} the space $H^0(V, K_V)$ has no trivial subrepresentation, there is no $G$-equivariant embedding  $T\hookrightarrow H^0(\OO_{\pp^7}(1))$.

Let $W$ be the cone over $\Si$ with vertex $P$. $W$ is a singular
threefold of degree 6  and it is mapped to itself by the
$G$-action. We now show that every nontrivial element of $G$ has
only isolated fixed points on $W$.

Notice that every fixed point of $g\in G$ on $W\setminus\{P\}$
lies on the line joining the vertex $P$ to a fixed point of $g$ on
$\Si$. 

Consider first an element $g\notin\Span{g_1}$. One checks
that $g$ has 3 isolated fixed points $Q_1, Q_2,Q_3$ on $\Si$, none
of which satisfies $z_{02}=0$. By assumption, $\chi_2(g)=\om^a$
with $a=1$ or $a=2$. Hence $g$ acts by multiplying the coordinate
$z_{01}$ by $\om^a$ and the coordinate $z_{02}$ by $\om^{2a}$.
Thus $g$ acts nontrivially on the line $PQ_i$, for $i=1,2,3$ and
$P,Q_1,Q_2,Q_3$ are the only fixed points of $g$ on $W$. Consider
now $g=g_1$. The fixed points of $g_1$ on $\Si$ are the six points
that lie on two of the $-1$-curves contained in $\Si$. If $e_0$ is
the strict transform of the line $x_0=0$ (say), then $e_0$ is
embedded in $\pp^6$ by the following map:
\begin{gather*}
z_{02}=0, \quad z_{12}= x_1, \quad z_{22}= x_2, \quad z_{10}=
\om^2 x_1,\\
z_{20}= \om^2 x_2, \quad z_{11}= \om x_1, \quad z_{21}=\om x_2,
\end{gather*}
The two fixed points of $g_1$ on $e_0$ correspond to the points
$Q_1$ of coordinates $x_1=1$, $x_2=0$ and $Q_2$ of coordinates
$x_1=0$, $x_2=1$. Hence the action of $g_1$ on the lines $PQ_i$ is
nontrivial (for $i=1,2$). Since the set of fixed points of $g_1$
on $\Si$ is the union of the $G$-orbits of $Q_1$ and $Q_2$, the
previous computation shows that the fixed points of $g_1$ on $W$
are $P$ and the points $Q_1,\dots,Q_6\in \Si$.

Consider now the subspace $T_1\subset H^0(\OO_W(3))$ on which the $G$-action  is trivial. If we set $z:=z_{01}$, then we can write:
\begin{equation}\label{dec}
H^0(\OO_W(3))=<z^3>\oplus z^2H^0(\OO_{\Si}(1))\oplus zH^0(\OO_{\Si}(2))\oplus H^0(\OO_{\Si}(3)),
\end{equation}
and this decomposition is compatible with the $G$-action. We identify $H^0(\OO_{\Si}(i))$ with the space  of homogeneous polynomials of degree $3i$ in $x_0, x_1,x_2$ that vanish of order at least $i$ at the coordinate points of $\pp^2$. Then the following is a basis of $T_1$:
\begin{gather*}
z^3, \ z(x_0^4x_1x_2+\omega^2x_1^4x_2x_0+\omega x_2^4x_0x_1), \  z(x_0^3x_1^3+\omega^2 x_1^3x_2^3 +\omega x_2^3x_0^3)\\
x_0^6x_1^3+x_1^6x_2^3+x_2^6x_0^3, \ x_0^6x_2^3+x_1^6x_0^3+x_2^6x_1^3, x_0^5x_1^2x_2^2+x_1^5x_2^2x_0^2+x_2^5x_0^2x_1^2\\
x_0^4x_1^4x_2+x_1^4x_2^4x_0+x_2^4x_0^4x_1, x_0^3x_1^3x_2^3.
\end{gather*}

 The dimension of $T_1$ is equal to 8  and the corresponding linear system $|T_1|$ is free,  since all the cubic powers of the homogeneous  coordinates $z_{ij}$ of $\pp^7$ restrict on $W$  to elements of $T_1$.

Let $V\in |T_1|$ be a surface with at
most rational double points and not passing through the points
fixed by a nontrivial element of $G$. Set $W_0:=\Proj(\OO_{\Si}
\oplus \OO_{\Si}(1))$ and let $H$ be the tautological bundle. The
system $|H|$ gives a morphism $W_0\to W$ that contracts the
divisor $E$ ``at infinity'' to the vertex $P$ of the cone and
restricts to an isomorphism $W_0\setminus E\to W\setminus\{P\}$.
The surface $V$ does not contain $P$, hence it can be identified
with an element of $|3H|$ on $W_0$ that we again denote by $V$.
Using adjunction on $W_0$, one sees that the canonical divisor
$K_V$ is the restriction of $H$ to $V$, hence $K^2_V=3H^3=18$. The
adjunction sequence gives: $p_g(V)=8$, $\chi(V)=9$, $q(V)=0$. Set
$S:=V/G$. Then by Proposition~\ref{campe}, $S$ is the canonical
model of a numerical Campedelli surface with $\pionealg=G$.

\begin{prop}\label{moduliB2} The family of surfaces of type
$B2$ depends on $6$ moduli.
\end {prop}

\begin{proof}
By the argument in the proof of Proposition \ref{moduliA}, it is enough to
 show that the  centralizer  $\Ga$ of $G$ in $\Aut W$  is discrete.

Let $\Ga_1$ be the connected component of the identity in $\Ga$ and let $\ga\in \Ga_1$.
The $\pp^6$ containing
$\Si$ is mapped to iself by $\ga$, and by continuity $\ga$ maps
each of the 6 lines of $\Si$ to itself. Hence the restriction of
$\ga$ to $\Si$ is induced by an element of $\Aut\pp^2$, that we
denote by $\ga_0$. Arguing as in the proof of
Proposition~\ref{moduliB1} one shows that $\ga_0$ is the identity
on $\pp^2$, hence $\ga$ restricts to the identity on
$\pp^6\supset\Si$. Thus the connected component of the identity in
$\Ga$ consists of the automorphisms that multiply $z_{01}$ by
$\lambda\in \C^*$ and do not change the remaining homogeneous
coordinates. Hence the surfaces of type B2 depend on $7-1=6$
moduli.
\end{proof}

\subsection{A common construction} \label{common}
We describe a common construction of the surfaces of type B1 and B2, suggested to us by Miles Reid.

Denote by $Z\subset \pp^8$ the Segre embedding of  $\pp^2\times{\pp^2}^*$ and let $G$ act on $Z\subset \pp^8$ as in the case of surfaces of type B1. It is easy to check that the corresponding representation on $H^0(\pp^8,\OO_{\pp^8}(1))$ is the sum of the 9 characters of $G$. 
We identify the group of characters  $G^*$ with $\Z_3^2$ as in \S\ref{ssecB2} and we denote by $z_{ij}$ homogeneous coordinates on $\pp^8$ such that  $G$ acts on $z_{ij}$ via the character $(i,j)$. The hyperplane $z_{00}=0$ cuts out on $Z$ the flag variety $x_0y_0+x_1y_1+x_2y_2=0$ and the hyperplane $z_{01}=0$ cuts out on $Z$ the invariant hypersurface $x_0y_0+\omega^2 x_1y_1+\om^2x_2y_2=0$. Set $\Si:=Z\cap\{z_{00}=z_{01}=0\}$. The surface $\Si$ is a smooth $G$-invariant Del Pezzo surface of degree 6 in $\pp^6$ and the $G$-action on $\Si$ coincides with the action given in \ref{Z3action}. 

Consider now homogeneous coordinates $z_{ij}, w$ on $\pp^9$, identify $\pp^8$ with the hyperplane $w=0$ and let $K$ be the cone over $Z$ with vertex $P:=[0,\dots 0, 1]$. Extend the $G$-action to $\pp^9$ by letting $G$ act on $w$ via the character $(0,1)$. 

Set $R_1:=H^0(\OO_{\pp^9}(3))^G$. The cubes $z_{ij}^3$, $w^3$ of the coordinates belong to $R_1$, hence the system $|R_1|$ is free.
Let now $V$ be the surface obtained by intersecting $K$ with  an element of $|R_1|$, the hyperplane $z_{00}=0$ and a hyperplane of the form $\lambda z_{01}+\mu w=0$.
The general $V$ is smooth by Bertini's theorem. Assume that $V$ has at most canonical singularities. If $\lambda=0$, then $V/G$ is a surface of type B2 and it is clear that every surface of type B2 can be obtained in this way. If $\lambda\ne 0$, then 
 $V/G$ is a surface of type B1 and every surface of type B1 can be obtained in this way, for instance taking $\mu=0$.

In particular, this construction proves the following:
\begin{prop}\label{limit}
Let $S$ be a surface of type B2. Then $S$ is a limit of surfaces of type B1.
\end{prop}
\begin{rem}{\em  In all three cases the general hyperplane section of $W$ is the smooth Del Pezzo surface $\Si$ of degree 6 in $\pp^6$, hence $W$ can be deformed to the cone over $\Si$. By the above construction,  in case B1 it is possible to preserve the $\Z_3^2$-action in the deformation. This is not possible in case A. Indeed the limit cone would have a smooth $\Z_9$-invariant hyperplane section, while $\pp^1\times\pp^1\times \pp^1$ has no such section.}
\end{rem}

\section{Geometry and moduli}
In this section we study some geometrical properties of numerical
Campedelli surfaces with $\pionealg$ of order 9 and of their
moduli space.
\subsection{Moduli}

Notice first of all the following consequence of Theorem~\ref{main}:
\begin{prop}\label{ample}
The canonical divisor $K_S$ of a general numerical Cam\-pe\-delli
surface $S$ with $\pionealg(S)$ of order $9$ is ample.
\end{prop}
\begin{proof}
The statement is immediate by Theorem~\ref{main} and by the
description of the three families of surfaces in \S3.
\end{proof}
The next result describes the moduli space of Campedelli surfaces with $|\pionealg|=9$:
\begin{thm}\label{moduli}
Let $\sM$ be the moduli space of numerical Campedelli surfaces,
let $\sM_A\subset \sM$ be the subset of surfaces with
$\pionealg=\Z_9$ and let $\sM_B\subset \sM$ be the subset of
surfaces with $\pionealg =\Z_3^2$. Then:
\begin{enumerate}

\item $\sM_A$ and $\sM_B$ are connected components of $\sM$; \item
$\sM_A$ is irreducible of dimension~$6$;

\item $\sM_B$ is irreducible of dimension~$7$.
\end{enumerate}
\end{thm}

\begin{proof}
Since $\pionealg$ is the profinite completion of $\pi_1$, it is a
topological invariant. Thus $\sM_A$ and $\sM_B$ are open and
closed in $\sM$. By Theorem~\ref{main},~(i), the points of $\sM_A$
correspond to the surfaces of type A, hence $\sM_A$ is irreducible
of dimension 6 by Proposition~\ref{moduliA}. By
Theorem~\ref{main},~(ii), the points of $\sM_B$ correspond to the
surfaces of type B1 and B2. By Proposition~\ref{limit}, the surfaces of type B2 lie in the
closure of the set of surfaces of type B1. In turn, the surfaces
of type B1 form an irreducible subset of dimension 7 by
Proposition~\ref{moduliB1}.
\end{proof}

\subsection{The topological fundamental group}
Having an explicit  construction of Campedelli surfaces with $|\pionealg|=9$ allows one to determine also the topological fundamental group.
\begin{prop}\label{pi1}
Let $S$ be a numerical Campedelli surface such that $\pionealg(S)$
has order 9\ {\rm (}$\Tors S$ has order 9\,{\rm )}. Then
$\pi_1(S)=\pionealg(S)$.
\end{prop}
\begin{proof}
The group $\Tors S$ is the largest abelian quotient of
$\pionealg(S)$. Thus, by Proposition~\ref{pi9}, the group $\Tors
S$ has order 9 if and only if $\pionealg(S)$ has order 9.

Assume that this is the case. Then it is enough to prove the
statement for one surface in each irreducible component of the
moduli space of numerical Campedelli surfaces with $\pionealg$ of
order 9. If $S$ is a surface of type A or B1 with $K_S$ ample,
then $V$ is a smooth ample divisor inside a smooth simply
connected threefold. Hence $V$ is simply connected by the
Lefschetz theorem on hyperplane sections and it is actually the
universal cover of $S$. By Theorem~\ref{moduli}, this proves the
statement. \end{proof}
\subsection{The bicanonical system}

We now study the bicanonical system.
\begin{thm}\label{basepoints}
Let $S$ be a numerical Campedelli surface such that $\pionealg(S)$
has order 9 and let $\Ga$ be the base locus of the bicanonical
system $|2K_S|$. Then:
\begin{enumerate}
\item if $S$ is of type A, then $\Ga$ consists of two points;
\item if $S$ is of type B1 then $\Ga$ is empty;
\item if $S$ is of type B2, then $\Ga$ consists of two points.
\end{enumerate}
\end{thm}

\begin{rem} \rm
Recall that, excluding the case $K^2=1$, $p_g=0$, the bicanonical
map of a minimal surface of general type is generically finite
{\rm(}\cite{xiaocan}{\rm)} and it is a morphism if either $p_g>0$
or $K^2>4$ (cf.\ \cite{ciro}{\rm)}. To our knowledge, surfaces of
type A and of type B2 are the only known examples of surfaces of
general type with $K^2>1$ whose bicanonical map is not a morphism.
\end{rem}

\begin{proof}[Proof of Theorem~\ref{basepoints}.] We use the
notation introduced in the previous sections.

Denote by $X$ the canonical model of $S$, denote by $p\colon V\to
X$ the quotient map and by $\epsi\colon S\to X$ the minimal
resolution. Since $|2K_S|=\epsi^*|2K_X|$, we study the base locus
of $|2K_X|$. Assume that $S$ is a surface of type A. In this case
$W$ is the Segre embedding of $\pp^1\times\pp^1\times\pp^1$ in
$\pp^7$ and a generator $g\in \Z_9$ acts on $W$ as \eqref{Z9action}. The $G$-action on $W$ induces an action on
$H^0(\OO_W(1))$ which is determined only up to multiplication by
an element of $G^*$. Let $\ze$ be a primitive $9$-th root of 1 and
denote by $\chi$ the character such that $\chi(g)=\zeta$. If we
require that the representation on $H^0(\OO_W(1))$ does not
contain the trivial character, then the representation is
determined uniquely and the following is a basis such that $G$
acts on $z_j$ as multiplication by $\chi^j$:
\begin{align*}
z_4:= z+\zeta y+\zeta^2x,&\quad z_7:=z+\zeta^4y+\zeta^8x,\\
z_1:=z+\zeta^7 y+\zeta^5x,&\quad
z_5:= xy+\zeta^7 yz+\zeta^8 xz, \\
z_8:=xy+\zeta^4 yz+\zeta^2xz,&\quad z_2:=xy +\zeta yz+\zeta^5xz,\\
z_3:=1,&\quad z_6:=xyz.
\end{align*}
By Proposition~\ref{Grep}, this choice of action on
$H^0(\OO_W(1))$ gives a $G$-isomor\-phism of $H^0(\OO_W(1))$ with
$H^0(K_V)=H^0(K_Y)$.

The threefold $W$ is projectively normal, hence $V$, being cut out on $W$ by a cubic hypersurface,  is also projectively normal. Thus $p^*H^0(2K_X)$ is generated by the restrictions to $V$ of
the four quadrics $z_iz_{9-i}$, for $i=1,\dots,4$. It is not
difficult to check that the zero locus on $W$ of these quadrics is
the union of the following curves:
\begin{gather}
L_1:=\{x=y=0\},\quad L_2:=\{ x=z=0\}, \quad L_3:=\{ y=z=0\}\\
L_4:=\{x=y=\infty\},\quad L_5:=\{ x=z=\infty\}, \quad L_6:=\{ y=z=\infty\} \nonumber
\end{gather}
The surface $V$ does not contain any of the curves $L_i$, since
the $G$ action on $V$ is free, while the points $(0,0,0)=L_1\cap
L_2\cap L_3$ and $(\infty, \infty,\infty)=L_4\cap L_5\cap L_6$ are
fixed points of $G$. Moreover, each of the $L_i$ is mapped to
itself by $g^3$. It follows that for every $i=1,\dots,6$ the set
$V\cap L_i$ consists of 3 distinct points, which form an orbit for
the action of $g^3$. In particular for $i=1,\dots,6$ the points of
$V\cap L_i$ are smooth for $V$. Hence $p^*|2K_X|$ has 18 base
points and $|2K_X|$ has two base points. Since the base points of
$|2K_X|$ are smooth points of $X$, also $|2K_S|$ has two base
points.

Assume now that $S$ is a surface of type B1 or B2, so that
$G=\Z_3^2$. In both cases the $G$-action on $W$ and on $\pp^7$ is
induced by the $G$-action on $\pp^2$ of \eqref{Z3action}. Let
$\chi_1$ be the character of $G$ such that $\chi_1(g_1)=\om$,
$\chi_1(g_2)=1$ and let $\chi_2$ be the character such that
$\chi_2(g_1)=1$, $\chi_2(g_2)=\om$. Denote by $(z_{ij})$, for
$(i,j)\in \Z_3^2\setminus\{0\}$, homogeneous coordinates on $\pp^7$ such that
$G$ acts on $z_{ij}$ as multiplication by $\chi_1^i\chi_2^j$.
Arguing as in the case of surfaces of type A, we have that
$p^*|2K_X|$ is generated by the restrictions to $V$ of the
following quadrics:
\begin{equation}\label{quadrics}
z_{10}z_{20},\enspace z_{01}z_{02},\enspace z_{11}z_{22},\enspace
z_{12}z_{21}.
\end{equation}
Assume that $S$ is of type B1. Then up to multiplication by
nonzero scalars we have the following equalities on $W$:
\begin{gather*}
z_{20}:=x_0y_0+\om x_1y_1+\om^2x_2y_2;\quad z_{10}:=x_0y_0+\om^2 x_1y_1+\om x_2y_2;\\
z_{02}:=x_0y_1+x_1y_2+x_2y_0;\quad z_{01}:=x_0y_2+x_1y_0+x_2y_1;\\
z_{22}:=x_0y_1+\om x_1y_2+\om^2x_2y_0;\quad z_{21}:=x_0y_2+\om x_1y_0+\om^2x_2y_1;\\
z_{12}:=x_0y_1+\om^2 x_1y_2+\om x_2y_0;\quad z_{11}:=x_0y_2+\om^2
x_1y_0+\om x_2y_1.
\end{gather*}
An easy computation shows that the zero locus $\Si_{ij}$ of
$z_{ij}$ on $W$ is a smooth surface for every $(0,0)\ne (i,j)\in
\Z_3^2$. Hence $\Si_{ij}$ is a smooth Del Pezzo surface of degree
6 and, in particular, it 
 contains six lines, that we denote by $e_i$ for
$i\in \Z_6$. The lines $e_i$ form an ``hexagon'', namely we can
arrange the indices in such a way that $e_i\cap e_j$ is a point if
$i-j=\pm 1$ and it is empty otherwise. The points $R_i:=e_i\cap
e_{i+1}$ are distinct for $i\in \Z_6$. The group $G$ preserves the
set $\{e_1,\dots,e_6\}$. Considering the intersection numbers
$e_ie_j$, one sees that an element of $G$ either maps each $e_i$
to itself or it induces a ``rotation'' of order 3 of the hexagon.
Notice that the whole group $G$ cannot act trivially on
$\{e_1,\dots,e_6\}$, since otherwise the points $R_1,\dots,R_6$
would be fixed by all the group.

Hence, if we denote by $\epsi_{ij}\colon \Si_{ij}\to\pp^2$ the birational
morphism that blows down $e_1, e_3, e_5$, then the $G$-action on
$\Si_{ij}$ descends to  a $G$-action on $\pp^2$. Since $G$ acts freely on $W$, and hence on $\Si_{ij}$, outside a finite set, by Lemma \ref{P2action} it follows that the action of $G$ on $\Si_{ij}$  is the one described in \S\ref{ssecB2}.
 In particular, the
$G$-invariant hyperplane sections of $\Si_{ij}$ are those defined
by the equations (\ref{map}). In particular all these sections are
irreducible except one, which is the union of the six lines of
$\Si_{ij}$. (It is not difficult to convince oneself that the
reducible section is given by $z_{2i\,2j}=0$).) So for every
choice of $(i_1,j_1),\dots, (i_4, j_4)\in \Z_3^2$ such that
$(i_r,j_r)+(i_s, j_s)\ne (0,0)$ for $1\le r,s\le 4$, the
intersection of $W$ with the subspace $z_{i_1j_1}=\cdots=
z_{i_4j_4}=0$ is a $G$-invariant set, properly contained in a
hyperplane section of the smooth elliptic curve $\Si_{i_1j_1}\cap
\Si_{i_2j_2}$. Hence this set has at most six points, all of which
have nontrivial stabilizer. So the base locus $\Ga_0$ on $W$ of
the quadrics (\ref{quadrics}) is a finite set (it is easy to check
that $\Ga_0$ is nonempty) and every point of $\Ga_0$ has
nontrivial stabilizer. Since $G$ acts freely on $V$, it follows
that $V\cap\Ga_0=\emptyset$. So the system $p^*|2K_X|$ is free
and, as a consequence, $|2K_X|$ and $|2K_S|$ are also free.

Assume now that $S$ is of type B2. We use the notation of \S\ref{ssecB2} and,
as in the case of surfaces of type B1, we study the zero set on
$W$ of the quadrics (\ref{quadrics}). Write $\Si$ for the
intersection of $W$ with $z_{01}=0$, namely $\Si$ is the only
$G$-invariant section of $W$ not containing the vertex $P$ of $W$.
Take $(i_1,j_1),\dots,(i_4, j_4)\in \Z_3^2$ such that
$(i_r,j_r)+(i_s, j_s)\ne (0,0)$ for $1\le r,s\le 4$. Up to
permuting the indices, we may assume that $(i_1, j_1)=(0,1)$ or
$(i_1, j_1)=(0,2)$. If $(i_1, j_1)=(0,1)$ then we may argue as in
the previous case and show that the intersection of $W$ with the
subspace $z_{01}=z_{i_2j_2}=z_{i_3j_3}=z_{i_4j_4}=0$ is a finite
set all of whose points have nontrivial stabilizer.

Now assume that $(i_1,j_1)=(0,2)$. Then the intersection of $W$
with the subspace $z_{02}=z_{i_2j_2}=z_{i_3j_3}=z_{i_4j_4}=0$ is
the join of $P$ and of the intersection of $\Si$ with the subspace
$z_{02}=z_{i_2j_2}=z_{i_3j_3}=z_{i_4j_4}=0$. By the formulae
(\ref{map}), as the indices $(i_2,j_2), (i_3,j_3), (i_4,j_4)$ vary
one obtains the six points of intersection of the six lines of
$\Si$. Summing up, the zero set on $W$ of the quadrics
(\ref{quadrics}) is the union of finitely many points with
nontrivial stabilizer and of six rulings of $W$. The same argument
that we have used for surfaces of type A shows that $|2K_X|$ has
two base points, that are smooth for $X$, and thus $|2K_S|$ also
has two base points.
 \end{proof}

\section{Proof of the classification}\label{secproof}

This section proves the classification theorem \ref{main}.
We use freely the notation and the assumptions of \S2. 

Recall that the universal cover   $Y$ of $S$  satisfies $K^2_Y=3p_g(Y)-6$. A 
 detailed study of surfaces satisfying this relation
 has been carried out by Konno (\cite{konnopg}).
We recall here some of his results. By Proposition~\ref{Y}, (ii),
the surface $Y$ belongs to type I in Konno's classification.
Denote by $V\subset \pp^7$ the image of the canonical map
$\fie\colon Y\to\pp^7$ and by $W$ the intersection of all the
quadrics of $\pp^7$ containing $V$. The natural linear $G$-action
on $H^0(Y,K_Y)$ descends to a $G$-action on
$\pp^7=\pp(H^0(Y,K_Y)^*)$ that preserves $V$ and $W$.

Following Fujita, we define the $\De$-genus of a projective
variety $W$ of $\pp^n$ as $\deg W+\dim W-n-1$.

\begin{prop}\label{delta01} \cite{konnopg}
The variety $W$ has dimension~$3$ and $\De$-genus $0$ or $1$.
Furthermore:
\begin{enumerate}
\item if $W$ has $\De$-genus $0$, then it is a rational normal
scroll and the ruling of $W$ induces a fibration $V\to \pp^1$
whose general fibre is a smooth plane quartic;

\item if $W$ has $\De$-genus $1$, then $W$ is normal and $V$ is
the intersection of $W$ with a cubic hypersurface. Moreover, in
this case the canonical ring of $Y$ is generated in degree~$1$,
hence the canonical map is a morphism and $V$ is the canonical
model of $Y$.
 \end{enumerate}
\end{prop}

\begin{proof} Cf. \cite{konnopg}, Theorem~3.1, Theorem~4.2,
Theorem~6.2. The statement on the normality of $W$ and on the
generation of the canonical ring of $Y$ is contained in the proof
of Theorem~3.1 of \cite{konnopg} (cf.\ 3.6, ibid.).
 \end{proof}

We will see that in our case $W$ has $\De$-genus $1$. In order to prove this we need the following:

\begin{lem}\label{pencil} The surface $V$ does not have a
$G$-invariant free pencil $|F|$ of curves of genus $g(F)\le 4$.
 \end{lem}

\begin{proof} Assume that such a pencil $|F|$ exists and denote
by $|F'|$ the pencil induced by $|F|$ on $S$.

Let $H$ be the subgroup of $G$ consisting of the elements that act
trivially on $|F|$ and let $h$ be the order of $H$. The general
$F'$ is isomorphic to $F/H$ for some $F$, hence
$g(F')=1+(g(F)-1)/h$ (recall that $H$ acts freely). Since
$g(F)\le4$, we either have $h=1$ and $g(F)=g(F')$ or $h=3$,
$g(F)=4$, $g(F')=2$.

If $h=1$, then $G=\Z_9$, since $\Aut |F| =\Aut\pp^1$ does not have
a subgroup isomorphic to $\Z_3^2$. But then $|F'|$ is a pencil of
genus $\le 4$ with two fibres of multiplicity $9$, contradicting
the adjunction formula.

If $h=3$, then the pencil $|F'|$ has two triple fibres,
corresponding to the two fixed points of the action of $G/H$ on
$\pp^1=|F|$. This again contradicts the adjunction formula, since
$g(F')=2$.
 \end{proof}

\begin{prop}\label{delta1}
The threefold $W\subset \pp^7$ has $\De$-genus $1$.
\end{prop}

\begin{proof} By Proposition~\ref{delta01} it is enough to
exclude that $W$ is a rational normal scroll. So assume by
contradiction that this is the case and denote by $f\colon
V\to\pp^1$ the fibration induced by the ruling of $W$. The general
fibre $F$ of $f$ is smooth of genus 3 by
Proposition~\ref{delta01}, (i), and by construction the $G$-action
on $V$ preserves the fibration $f$. This contradicts
Lemma~\ref{pencil}.
\end{proof}

Since $W$ has $\De$-genus 1 by Proposition~\ref{delta1}, it
follows by Proposition~\ref{delta01}, (ii) that the canonical
image $V$ of $Y$ is the intersection of $W$ with a cubic
hypersurface. As a consequence, we get the following:
\begin{cor}\label{fix} The $G$-action on $\pp^7$ restricts to
an action on $W$ which is free outside a finite subset of $W$.
\end{cor}

\begin{proof} As already remarked, $G$ acts on
$V\subset W\subset\pp^7$ compatibly. By Proposition~\ref{delta01}
and Proposition~\ref{delta1}, the surface $V$ is the canonical
model of $Y$ and thus $G$ acts freely on $V$. Since by
Proposition~\ref{delta01} (ii) the divisor $V\subset W$ is very
ample, it follows that the set of points of $W$ with nontrivial
stabilizer has dimension $\le 0$.
\end{proof}

We also need the following two elementary results.

 \begin{lem}\label{P2action}
Assume that $G=\Z_3^2$ acts on $\pp^2$ such that the action is
free outside a finite set. Then we can choose generators
$g_1,g_2\in G$ and homogeneous coordinates $(x_0,x_1,x_2)$ on
$\pp^2$ so that the action is as follows:
\[
g_1\colon(x_0,x_1,x_2)\mapsto(x_0,\om x_1, \om^2x_2), \quad
g_2\colon(x_0,x_1,x_2)\mapsto(x_1,x_2,x_0),\quad
\]
where $\om\ne1$ is a third root of 1. \end{lem}

\begin{proof} The fixed locus of a nontrivial element
$g_1\in G$ consists of three points. In suitable homogeneous
coordinates the action of $g_1$ is given by:
\[
g_1\colon (x_0,x_1,x_2)\mapsto(x_0,\om x_1,\om^2x_2).
\]
Now consider $g_2\in G\setminus\Span{g_1}$. The element $g_2$ acts
on the fixed points of $g_1$. If this action is trivial, then
$g_2$ acts by:
\[
g_2\colon(x_0,x_1,x_2)\mapsto(x_0,\om^a x_1, \om^bx_2),
\]
for some $a,b\in \Z_3\setminus\{0\}$. Then either $g_1g_2$ or
$g_1g_2^2$ fixes a line pointwise, contradicting the assumptions.
Hence we conclude that $g_2$ permutes the fixed points of $g_1$
cyclically. Up to rescaling the coordinates, and possibly
replacing $g_2$ by $g_2^2$, the action can be written as stated.
\end{proof}

 \begin{lem}\label{P1action}
Let $G$ be a group of order~$9$ that acts on
$\pp^1\times\pp^1\times \pp^1$ such that the action is free
outside a finite subset. If $G$ permutes the three factors of
$\pp^1\times\pp^1\times\pp^1$ in a nontrivial way, then $G=\Z_9$
and there are affine coordinates $x,y,z$ on the three copies of
$\pp^1$ such that a generator $g$ of\/ $G$ acts by:
\[
g\colon(x,y,z)\mapsto(y,z,\om x),
\]
where $\om$ is a primitive third root of $1$.
\end{lem}

\begin{proof} The $G$-action permutes the three factors of
$\pp^1\times\pp^1\times \pp^1$, hence induces a homomorphism
$\psi\colon G\to S_3$, which is nontrivial by assumption. Thus the
image of $\psi$ has order~$3$.

Assume that $G=\Z_3^2$ and let $g\in G$ be such that
$\psi(g)=(123)$. Then there are affine coordinates $x.y,z$ on the
three copies of $\pp^1$ such that $g$ acts as follows:
\[
g\colon(x,y,z)\mapsto(y,z,x).
\]
Hence $g$ fixes pointwise the diagonal $\{(P,P,P):P\in \pp^1\}$,
contradicting the assumptions.

Thus $G=\Z_9$. Let $g\in G$ be an element such that
$\psi(g)=(123)$ (notice that $g$ generates $G$). Then in suitable
coordinates the action of $g$ can be written as claimed.
\end{proof}

The
$G$-action on $H^0(Y,K_Y)=H^0(\pp^7,
\OO_{\pp^7}(1))=H^0(W,\OO_W(1))$ induces a linearization of the
line bundle $\OO_W(3)$ and therefore a decomposition
\[
H^0(W,\OO_W(3))=\bigoplus_{\chi\in G^*}T_{\chi},
\]
where $T_{\chi}$ denotes the eigenspace corresponding to the
character $\chi$. In particular, $T_1$ is the subspace of
$G$-invariant vectors. The following remark will be useful in
determining the subsystem of $|\OO_W(3)|$ which parametrizes the
surfaces $V$.
\begin{lem}\label{system}
Let $P\in W$ be a point which is fixed by some nontrivial element
$g\in G$. If the system $|T_1|$ is free, then $P$ is a base point
of $|T_{\chi}|$ for every $\chi$ such that $\chi(g)\ne 1$.
\end{lem}
\begin{proof} Let $\chi\in G^*$ be a character such that
$\chi(g)\ne 1$. Fix a section $\si_0\in T_1$ such that
$\si_0(P)\ne 0$. For any $\si\in |T_{\chi}|$ consider the rational
function $f_{\si}:=\frac{\si}{\si_0}$. The function $f_{\si}$ is
defined at $P$ and one has:
\[
f_{\si}(P)=f_{\si}(g(P))=\chi(g\inv)f_{\si}(P).
\]
It follows that $f_{\si}(P)=0$ and thus also $\si(P)=0$.
\end{proof}

Varieties of $\De$-genus 1 have been classified by Fujita (cf.\
\cite{fujita1}, \cite{fujita2}, \cite{fujita3}). We recall that by
Proposition~\ref{delta01}, (ii) the threefold $W$ is normal. Here
is the list of normal threefolds of $\De$-genus 1 in $\pp^7$ (cf.\
\cite{fujita3}):
\begin{itemize}

\item[1)] the cone over a (weak) Del Pezzo surface $\Si\subset
\pp^6$ of degree 6;

\item[2)] the Segre embedding of $\pp^1\times\pp^1\times \pp^1$
into $\pp^7$;

\item[3)] the image of the flag variety $\{x_0y_0+x_1y_1+x_2y_2=0\}\subset \pp^2\times{\pp^2}^*$ under the Segre embedding of $\pp^2\times{\pp^2}^*$.

\item[4)] four singular examples that are not cones. In the notation of \cite{fujita3}, these are the cases (vi), (si31i), (si22i),  (si211).
\end{itemize}
To prove Theorem~\ref{main} we examine cases 1)--4) separately. We
start by showing that case 4) does not occur.

\medskip\noindent{\bf Case 4)} Consider first  case 4). If $W$ is of type (vi), (si31i) or  (si22i), then  by \cite[Theorem~2.9]{fujita3} the singular locus of $W$
consists of a line $r$ and, in case (vi), possibly also of an
isolated double point. Hence the line $r$ is mapped to itself by
the action of $G$ on $\pp^7$. Then $Z:=V\cap r$ is the
intersection of $r$ with a cubic hypersurface and $G$ acts freely
on $Z$. This is not possible, since $Z$ is either equal to $r$ or
it consists of at most 3 points.

If $W$ is of type (si211), then by \cite{fujita3} it is contained in the cone
$K$ over $M:=\pp(\OO_{\pp^1}(1)\oplus\OO_{\pp^1}(1)
\oplus\OO_{\pp^1}(2))$ embedded in $\pp^6$ by its tautological
system. We denote by $P$ the vertex of $K$. By
\cite[Theorem~2.9]{fujita3} $P$ is the only isolated double point
of $W$, hence $P$ is a fixed point of $G$. Notice, in particular,
that $V$ does not contain $P$. Let $\pp^6\subset \pp^7$ be the
$G$-invariant hyperplane not containing $P$ and let $q\colon
\pp^7\to\pp^6$ be the projection with center $P$. The threefold
$q(W)$ is $G$-invariant and it is isomorphic to $M$. Hence the map
$V\to\pp^1$ obtained by composing the projection $q$ with the
ruling $M\to\pp^1$ gives a free pencil of curves $|F|$ which is
acted on by $G$. A general $F$ is the intersection in $\pp^3$ of a
cubic (given by the cubic equation of $V$) and a quadric through
$P$ (corresponding to $W$) and therefore it has genus 4. Hence we
have a contradiction to Lemma~\ref{pencil}.

\medskip\noindent{\bf Case 1)} We show that this case corresponds
to surfaces of type B2.

First of all we show that we have $G=\Z_3^2$ in this case. The
vertex $P$ of the cone $W$ is a fixed point of the $G$-action on
$\pp^7$. In particular, $V$ does not contain $P$.

Since the $G$-action on $\pp^7$ is induced by a linear
representation of $G$ on $H^0(\OO_{\pp^7}(1))=H^0(Y,K_Y)$, there
exists a $G$-invariant $\pp^6\subset \pp^7$ that does not contain
$P$. Let $\Si:=W\cap \pp^6$. Then $\Si$ is a (weak) Del Pezzo
surface of degree 6 preserved by the $G$-action.

The projection from $\pp^7\to \pp^6$ with center $P$ restricts to
a $G$-invariant degree three morphism $p\colon V\to \Si$. If a
point $Q\in \Si$ is fixed by every element of $G$, then $G$
permutes the points in the fibre $p\inv(Q)$, which consists of at
most three points. Since $G$ has order 9, this contradicts the
assumption that $G$ acts freely on $V$. In particular, $G$ cannot
be cyclic, since every automorphism of a rational surface has at
least one fixed point. So $G=\Z_3^2$.

The same argument shows that the Del Pezzo surface $\Si$ is
smooth. Indeed, by \cite[Theorem~8]{nagata}, $\Si$ is the blowup
of $\pp^2$ along a curvilinear scheme $Z$ of dimension $0$ and
length 3. If $Z$ is not reduced, then $\Si$ is singular and it has
precisely one singular point, which is necessarily fixed by every
element of $G$.

Arguing as in the proof of Theorem \ref{basepoints}, one shows that there is a birational morphism $\epsi\colon\Si\to\pp^2$ such that $\epsi$ contracts three disjoint lines of $\Si$ to the coordinate points of $\pp^2$ and the $G$-action on $\Si$ descends to a $G$-action on $\pp^2$.  Assume that this last action is
not free outside a codimension 2 subset of $\pp^2$. Then there
exists $g\in G$, $g\neq 0$, that fixes a line $r$ of $\pp^2$
pointwise. Since $G$ is abelian, we have $g'(r)=r$ for every
$g'\in G$. In particular, $r$ does not contain any of the
exceptional points of $\epsi\inv$, since these are linearly
independent and they are permuted cyclically by $G$. Now let
$g'\in G$ be such that $g,g'$ generate $G$ and let $A\in r$ be a
point such that $g'(A)=A$. Then the point $A$ is fixed by all the
elements of $G$ and therefore the point $\epsi\inv(A)\in \Si$ is
also fixed by every element of $G$, contradicting the remarks
above.

Thus we conclude that $G$ acts freely on $\pp^2$ outside a finite
set and that for  a suitable choice of homogeneous coordinates of
$\pp^2$ the $G$-action can be written as in Lemma~\ref{P2action}.

Hence we may assume that the $G$-action on $\pp^7$ is the one
described in \S\ref{ssecB2}. The point $P$, being invariant for all the group
$G$, is a coordinate point of $\pp^7$. Let $z_{ij}$ denote the
only coordinate that does not vanish at $P$, so that the invariant
$\pp^6\subset \pp^7$ that does not contain $P$ is defined by
$z_{ij}=0$. The representation of $G$ on $H^0(\OO_{\pp^6}(1))$
contains all the nontrivial characters of $G$ except
$\chi_1^i\chi_2^j$ (see \S\ref{ssecB2} for the notation). In turn,
$H^0(\OO_{\pp^6}(1))$ is isomorphic to the subspace $T\subset
H^0(\OO_{\pp^2}(3))$ of cubics vanishing at the coordinate points
of $\pp^2$ and the $G$-action on $T$ induced by this isomorphism
differs from the action given in (\ref{base}) by multiplication by
a character of $G$. Hence we have $(i,j)=(0,1)$ or $(i,j)=(0,2)$.
It follows that, possibly up to replacing $g_1$ by $g_1^2$, the
embedding $\Si\subset \pp^6\subset \pp^7$ is induced by the
rational map (\ref{map}). By Proposition~\ref{delta01}, (ii), $V$
is an element of $|\OO_{W}(3)|$ that is $G$-invariant, hence there
is a $\chi\in G^*$ such that $V\in |T_{\chi}|$, where
$T_{\chi}\subset H^0(\OO_W(3))$ is the eigenspace corresponding to
$\chi$. We have seen in \S\ref{ssecB2} that the system $|T_1|$ is free.
Hence by Lemma~\ref{system} the vertex $P$ of $W$ is in the base
locus of the system $|T_{\chi}|$ for every $1\ne \chi\in G^*$.
Since $G$ acts freely on $V$, it follows that $V$ belongs to
$|T_1|$ and the minimal desingularization of $S:=V/G$ is a surface
of type B2.

\medskip\noindent{\bf Case 2)} We show that this case
corresponds to surfaces of type A. We use the notation of \S3 for
the homogeneous and affine coordinates on
$W=\pp^1\times\pp^1\times\pp^1$.

For $i=1,2,3$, the projection onto the $i$-th factor $p_i\colon
W\to \pp^1$ restricts on $V$ to a free pencil of genus 4. By
Lemma~\ref{pencil}, these pencils are not $G$-invariant, hence
there is at least an element of $G$ that permutes them. Since $G$
acts freely outside a finite subset, by Lemma~\ref{P1action} we
have $G=\Z_9$ and $G$ acts as in Lemma~\ref{P1action}.

By Proposition~\ref{delta01}, (ii), $V$ is an element of
$|\OO_{W}(3)|$ that is $G$-invariant, hence there is a $\chi\in
G^*$ such that $V\in |T_{\chi}|$, where $T_{\chi}\subset
H^0(\OO_W(3))$ is the eigenspace corresponding to $\chi$. We have
seen in \S\ref{ssecA} that the system $|T_1|$ is free. Let $Q\in W$ be the
point with affine coordinates $x=y=z=0$. The point $Q$ is fixed by
all the group $G$, hence, by Lemma~\ref{system}, $Q$ is in the
base locus of the system $|T_{\chi}|$ for every $1\ne \chi\in
G^*$. Since $G$ acts freely on $V$, it follows that $V$ belongs to
$|T_1|$ and the minimal desingularization of $S:=V/G$ is a surface
of type A.

\medskip\noindent{\bf Case 3)} We show that this case corresponds to
surfaces of type B1.

Here $W$ is the flag variety
$\{x_0y_0+x_1y_1+x_2y_2=0\}\subset \pp^2\times{\pp^2}^*$, embedded
in $\pp^7$ as a hyperplane section of the Segre embedding
$\pp^2\times{\pp^2}^*\into\pp^8$ (cf. \S \ref{ssecB1}).
By Lemma \ref{autW}, the action of $G$ on $W$ is induced by a $G$-action on $\pp^2$.
The fixed points on $W$ of an element $g\in \Aut\pp^2$ correspond
to pairs $(P,r)$ where $P\in \pp^2$ is a fixed point of $g$, $r\in
{\pp^2}^*$ is a fixed line and $P\in r$. Since $G$ acts freely
outside a finite subset of $W$, the elements of $G$ have finitely
many fixed points on $\pp^2\times{\pp^2}^*$. Hence the $G$-action
on $\pp^2$ is free outside a finite subset of $\pp^2$. Assume that $G=\Z_9$ and let $g\in G$ be a generator. Then there
are homogeneous coordinates such that $g$ acts by:
\[
(x_0,x_1,x_2)\mapsto (x_0, \zeta x_1, \zeta ^k x_2),
\]
where $\zeta$ is a primitive $9$-th root of 1. Then $x_0y_0$ and
$x_1y_1$ restrict on $W$ to independent sections of $\OO_W(1,1)$
that belong to the same $G$-eigenspace. This contradicts
Proposition~\ref{Grep}. Hence $G=\Z_3^2$ and the action of $G$ on
$\pp^2$ can be written as in Lemma~\ref{P1action}.

By Proposition~\ref{delta01}, (ii), $V$ is an element of
$|\OO_{W}(3,3)|$ that is $G$-invariant, hence there is a $\chi\in
G^*$ such that $V\in |T_{\chi}|$, where $T_{\chi}\subset
H^0(\OO_W(3,3))$ is the eigenspace corresponding to $\chi$. We
have seen in \S\ref{ssecB1} that the system $|T_1|$ is free. We observe that
for every $0\ne g\in G$ there is a point $Q\in W$ such that
$gQ=Q$. Thus by Lemma~\ref{system}, $Q$ is in the base locus of
the system $|T_{\chi}|$ if $\chi(g)\ne 1$. Since $G$ acts freely
on $V$, it follows that $V$ belongs to $|T_1|$ and the minimal
desingularization of $S:=V/G$ is a surface of type B1.

\bigskip

\bigskip\noindent
\begin{minipage}{12.5cm}
\parbox[t]{5.7cm}{Margarida Mendes Lopes\\
Departamento de Matem\'atica\\
Instituto Superior T\'ecnico\\
Universidade T{\'e}cnica de Lisboa\\
Av.~Rovisco Pais\\
1049-001 Lisboa, PORTUGAL\\
mmlopes@math.ist.utl.pt
 } \hfill
\parbox[t]{5.5cm}{Rita Pardini\\
Dipartimento di Matematica\\
Universit\`a di Pisa\\
Largo B. Pontecorvo, 5\\
56127 Pisa, Italy\\
pardini@dm.unipi.it}
\end{minipage}

\end{document}